# Uniform error bounds for smoothing splines

## P. P. B. Eggermont[1] and V. N. LaRiccia[1]

*University of Delaware*

**Abstract:** Almost sure bounds are established on the uniform error of smoothing spline estimators in nonparametric regression with random designs. Some results of Einmahl and Mason (2005) are used to derive uniform error bounds for the approximation of the spline smoother by an "equivalent" reproducing kernel regression estimator, as well as for proving uniform error bounds on the reproducing kernel regression estimator itself, uniformly in the smoothing parameter over a wide range. This admits data-driven choices of the smoothing parameter.

## 1. Introduction

In this paper, we study uniform error bounds for the smoothing spline estimator of arbitrary order for a nonparametric regression problem. In effect, we approximate the smoothing spline by a kernel-like estimator, and give sharp bounds on the approximation error under very mild conditions on the nonparametric regression problem, as well as on the uniform error on the kernel-like estimator. An application to obtaining confidence bands is pointed out.

Let $(X_1, Y_1), (X_2, Y_2), \ldots, (X_n, Y_n)$ be a random sample of the bivariate random variable $(X, Y)$ with $X \in [0, 1]$, almost surely. Assume that

$$(1.1) \qquad f_o(x) = \mathbb{E}[Y \mid X = x]$$

exists, and that for some natural number $m$,

$$(1.2) \qquad f_o \in W^{m,\infty}(0, 1),$$

where for $a < b$, the Sobolev spaces $W^{m,p}(a, b)$, $1 \le p \le \infty$, are defined as

$$(1.3) \qquad W^{m,p}(a, b) = \left\{ f \in C^{m-1}[a, b] \;\middle|\; \begin{array}{l} f^{(m-1)} \text{ abs. continuous} \\ f^{(m)} \in L^p(a, b) \end{array} \right\},$$

see, e.g., [2].

Regarding the design, assume that

$$(1.4) \qquad \begin{array}{l} X_1, X_2, \ldots, X_n \text{ are independent and identically distributed,} \\ \text{having a probability density function } w \text{ with respect to} \\ \text{Lebesgue measure on } (0, 1), \end{array}$$

[1]Food and Resource Economics, University of Delaware, Newark, Delaware 19716-1303, e-mail: eggermon@udel.edu; lariccia@udel.edu

and that

$$(1.5) \qquad w_1 \leq w(t) \leq w_2 \quad \text{for all} \quad t \in [\, 0\,,\, 1\,],$$

for positive constants $w_1$ and $w_2$.

With the random variable $(X, Y)$, associate the noise $D$ by

$$(1.6) \qquad D = Y - f_o(X),$$

and define $D_i = Y_i - f_o(X_i), \; i = 1, 2, \ldots, n$. Assume that

$$(1.7) \qquad \sup_{x \in [\, 0\,,\, 1\,]} \mathbb{E}[\,|\, D\,|^{\kappa} \mid X = x\,] < \infty \quad \text{for some } \kappa > 2.$$

(With the assumption (1.2), this is equivalent to $\sup_x \mathbb{E}[\,|\, Y\,|^{\kappa} \mid X = x\,] < \infty$.)

Under the above conditions, uniform error bounds for the Nadaraya-Watson estimator have been established by Deheuvels and Mason [8] for a random choice of the smoothing parameter, and by Einmahl and Mason [12] uniformly in the smoothing parameter over a wide range. We recall that the Nadaraya-Watson estimator is defined as

$$\widehat{f_n}(t) = \tfrac{1}{n} \sum_{i=1}^{n} Y_i \, K_h(t - X_i) \Big/ \tfrac{1}{n} \sum_{i=1}^{n} K_h(t - X_i),$$

where, $K_h(t) = h^{-1} K(h^{-1} t)$ for some nice "kernel" $K$. In this case, $\widehat{f_n(x)}$ is an estimator of $f_o(x) = \mathbb{E}[\, Y \mid X = x\,]$. For some earlier results on uniform error bounds for Nadaraya-Watson estimators, see, e.g., [16] and [15].

For the smoothing spline estimator, we must come to terms with the fact that the estimator is defined implicitly as the solution, denoted by $f = f^{nh}$, of a minimization problem,

$$(1.8) \qquad \begin{aligned} &\text{minimize} \quad \mathrm{IS}(f) \overset{\text{def}}{=} \tfrac{1}{n} \sum_{i=1}^{n} |\, f(X_i) - Y_i\,|^2 + h^{2m} \,\|\, f^{(m)}\,\|^2 \\ &\text{subject to} \quad f \in W^{m,2}(\, 0\,,\, 1\,), \end{aligned}$$

where $\|\cdot\|$ denotes the $L^2(0,1)$ norm. Thus, the bulk of the paper is devoted to establishing that for all $t \in [\, 0\,,\, 1\,]$,

$$(1.9) \qquad f^{nh}(t) - \mathbb{E}[\, f^{nh}(t) \mid X_1, \ldots, X_n\,] = \tfrac{1}{n} \sum_{i=1}^{n} D_i \, \mathfrak{R}_{wmh}(X_i, t) + \varepsilon^{nh}(t),$$

where $\|\, \varepsilon^{nh}\,\|_{\infty}$ is negligible compared to the leading term in (1.9), and $\mathfrak{R}_{wmh}$ is the Green's function for a suitable boundary value problem, see (2.13). Here, $\|\cdot\|_{\infty}$ denotes the $L^{\infty}(\, 0\,,\, 1\,)$ norm. The approach taken follows Eggermont and LaRiccia [10].

The precise results are as follows. For $\gamma > 0$, define the intervals

$$(1.10) \qquad \mathcal{H}_n(\gamma) = \Big[\, \gamma \,\big(\tfrac{\log n}{n}\big)^{1-2/\kappa},\, \tfrac{1}{2}\,\Big],\; \mathcal{G}_n(\gamma) = \Big[\, \gamma \,\big(\tfrac{\log n}{n}\big)^{1-2/\lambda},\, \tfrac{1}{2}\,\Big],$$

where $\lambda$ is unspecified but satisfies $2 < \lambda < \min(\kappa, 4)$.

**Theorem 1.** *Under the assumptions* (1.4) - (1.7) *on the model* (1.1), *the error term* $\varepsilon^{nh}$ *in* (1.9) *satisfies almost surely,*

$$T_{UE}(\gamma) \stackrel{\text{def}}{=} \limsup_{n \to \infty} \sup_{h \in \mathcal{H}_n(\gamma)} \frac{\|\varepsilon^{nh}\|_\infty}{h^{-1/2}(nh)^{-1}\{\log(1/h) \vee \log\log n\}} < \infty.$$

The uniform-in-bandwidth character of this theorem (which admits *random* choices of the smoothing parameter) stands out. Regarding the actual error bound, if $h \in \mathcal{G}_n(\gamma)$, then $h \gg (n^{-1}\log n)^{1/2}$ and the error term in (1.9) can be ignored. Note that for $m \geq 2$ and $\kappa > 3$, this covers the optimal $h$, which behaves like $(n^{-1}\log n)^{1/(2m+1)}$. The theorem makes the smoothing spline more accessible as an object of study. Here, we consider uniform error bounds on the estimator. For *cubic* smoothing splines in a somewhat different setting, uniform error bounds were derived by Chiang, Rice and Wu [5].

**Main Theorem.** *Assume the conditions* (1.2) *through* (1.7) *on the model* (1.1). *Then, the spline estimator of order* $m$ *satisfies almost surely,*

$$Q_{UE}(\gamma) \stackrel{\text{def}}{=} \limsup_{n \to \infty} \sup_{h \in \mathcal{G}_n(\gamma)} \frac{\|f^{nh} - f_o\|_\infty}{\sqrt{h^{2m} + (nh)^{-1}\{\log(1/h) \vee \log\log n\}}} < \infty.$$

The constant $Q_{UE}$ depends on the unknown regression function $f_o$ through the bias. If we restrict $h$ such that $h \ll (n^{-1}\log n)^{1/(2m+1)}$, then this dependence disappears, e.g., if for $m \geq 2$ and $\kappa > 2 + (1/m)$, we let

(1.11)                    $$\mathcal{F}_n(\gamma) = \left[\gamma\left(n^{-1}\log n\right)^{1-2/\kappa}, \ n^{-1/(2m+1)}\right],$$

then

(1.12)        $$\mathfrak{Q}_{UE}(\gamma) \stackrel{\text{def}}{=} \limsup_{n \to \infty} \sup_{h \in \mathcal{F}_n(\gamma)} \frac{\|f^{nh} - f_o\|_\infty}{\sqrt{(nh)^{-1}\{\log(1/h) \vee \log\log n\}}} < \infty,$$

and $\mathfrak{Q}_{UE}$ does not depend on $f_o$. This has obvious consequences for the construction of confidence bands. Since it seems reasonable that the value of $\mathfrak{Q}_{UE}$ can be determined via bootstrap techniques, then almost sure confidence bands in the spirit of Deheuvels and Mason [8] and the Bonferroni bounds of Eubank and Speckman [14] may be obtained. The full import of this will be explored elsewhere.

## 2. The smoothing spline estimator

Let $m \in \mathbb{N}$ and $h > 0$ be fixed. The smoothing spline estimator, denoted by $f^{nh}$, is defined as the solution of the minimization problem (1.8). The problem (1.8) always has solutions, and for $n \geq m$, the solution is unique, almost surely. For more on spline smoothing, see, e.g., [13] or [24].

A closer look at the spline smoothing problem reveals that $f(X_i)$ is well-defined for any $f \in W^{m,2}(0,1)$. In particular, there exists a constant $c$ such that for all $f \in W^{m,2}(0,1)$ and all $x \in [0,1]$,

(2.1)                        $$|f(x)| \leq c\left\{\|f\|^2 + \|f^{(m)}\|^2\right\}^{1/2},$$

see, e.g., [2]. Then, a simple scaling argument shows that there exists a constant $c_m$ such that for all $0 < h \le 1$, all $f \in W^{m,2}(0,1)$, and all $t \in [0,1]$,

$$(2.2) \qquad |f(t)| \le c_m \, h^{-1/2} \, \|f\|_{mh}^2.$$

Here,

$$(2.3) \qquad \|f\|_{mh} \overset{\text{def}}{=} \big\{ \, \|f\|^2 + h^{2m} \, \|f^{(m)}\|^2 \, \big\}^{1/2}.$$

Of course, the inequality (2.2) is geared towards the uniform design. For the present, "arbitrary" design, it is more appropriate to consider the inner products

$$(2.4) \qquad \langle \, f, g \, \rangle_{wmh} = \langle \, f, g \, \rangle_{L^2((0,1),w)} + h^{2m} \, \langle \, f^{(m)}, g^{(m)} \, \rangle_{L^2(0,1)},$$

where $\langle \, \cdot \, , \, \cdot \, \rangle_{L^2(0,1)}$ is the usual $L^2(0,1)$ inner product and

$$(2.5) \qquad \langle \, f, g \, \rangle_{L^2((0,1),w)} = \int_0^1 f(t) \, g(t) \, w(t) \, dt.$$

The norms are then defined by $\|f\|_{wmh} = \big\{ \, \langle \, f, f \, \rangle_{wmh} \, \big\}^{1/2}$. With the design density being bounded and bounded away from zero, see (1.5), it is obvious that the norms $\| \cdot \|_{mh}$ and $\| \cdot \|_{wmh}$ are equivalent, uniformly in $h$. In particular, with the constants $w_1$ and $w_2$ as in (1.5), for all $f \in W^{m,2}(0,1)$,

$$(2.6) \qquad w_1 \, \|f\|_{mh} \le \|f\|_{wmh} \le w_2 \, \|f\|_{mh}.$$

(Note that, actually, $w_1 \le 1 \le w_2$.) Then, the analogue of (2.3) holds: There exists a constant $c_m$ such that for all $0 < h \le 1$, all $f \in W^{m,2}(0,1)$, and all $t \in [0,1]$,

$$(2.7) \qquad |f(t)| \le c_m \, h^{-1/2} \, \|f\|_{wmh}.$$

For later use, we quote the following multiplication result which follows readily with Cauchy-Schwarz : There exists a constant $c$ such that for all $f$ and $g \in W^{1,2}(0,1)$,

$$(2.8) \qquad \|f \, g\|_{L^1((0,1),w)} + h \, \|(f \, g)'\|_{L^1(0,1)} \le c \, \|f\|_{w,1,h} \, \|g\|_{w,1,h}.$$

Also, there exist constants $c_{k,k+1}$ such that for all $f \in W^{k+1,2}(0,1)$,

$$(2.9) \qquad \|f\|_{w,k,h} \le c_{k,k+1} \, \|f\|_{w,k+1,h}.$$

The inequality (2.7) says that the linear functionals $f \mapsto f(t)$ are continuous in the $\| \cdot \|_{wmh}$-topology, so that $W^{m,2}(0,1)$ with the inner product $\langle \, \cdot \, , \, \cdot \, \rangle_{wmh}$ is a reproducing kernel Hilbert space, see [3]. Thus, by the Riesz-Fischer theorem on the representation of bounded linear functionals on Hilbert space, for each $t$, there exists an element $\mathfrak{R}_{wmht} \in W^{m,2}(0,1)$ such that for all $f \in W^{m,2}(0,1)$,

$$(2.10) \qquad f(t) = \langle \, f, \mathfrak{R}_{wmht} \, \rangle_{wmh}.$$

Applying this to $\mathfrak{R}_{wmht}$ itself gives $\mathfrak{R}_{wmht}(s) = \langle \, \mathfrak{R}_{wmht}, \mathfrak{R}_{wmhs} \, \rangle_{wmh}$, so that it makes sense to define

$$(2.11) \qquad \mathfrak{R}_{wmh}(t,s) = \mathfrak{R}_{wmht}(s) = \mathfrak{R}_{wmhs}(t) \quad \text{for all } s, t \in [0,1].$$

Then, again the inequality (2.7) implies that

$$(2.12) \qquad \| \, \mathfrak{R}_{wmh}(\, t \, , \, \cdot \, ) \, \|_{wmh} \leq c_m \, h^{-1/2},$$

with the same constant $c_m$.

Finally, we observe that reproducing kernels may be interpreted as the Green's functions for appropriate boundary value problems, see, e.g., [9]. In the present case, $\mathfrak{R}_{wmh}(\, t \, , s)$ is the Green's function for

$$(2.13) \qquad \begin{aligned} (-h^2)^m \, u^{(2m)} + w \, u = v \quad &\text{on} \quad (\, 0 \, , 1 \, ), \\ u^{(k)}(0) = u^{(k)}(1) = 0, \quad &k = m, \ldots, 2m - 1. \end{aligned}$$

In case $w(\, t \, ) = 1$ for all $t$ (the uniform density), we denote $\mathfrak{R}_{wmh}$ by $\mathcal{R}_{mh}$.

We finish this section by showing that the little information we have on the reproducing kernels suffices to prove some useful bounds on random sums of the form

$$\tfrac{1}{n} \sum_{j=1}^{m} D_i \, f(X_i),$$

with $D_1, D_2, \ldots, D_n$ and $X_1, X_2, \ldots, X_n$ as in Section 1, and $f \in W^{m,2}(\, 0 \, , 1 \, )$ random, i.e., depending on the $D_i$ and $X_i$. To obtain these bounds, let

$$(2.14) \qquad \mathfrak{S}_{nh}(t) \stackrel{\text{def}}{=} \tfrac{1}{n} \sum_{i=1}^{n} D_i \, \mathfrak{R}_{wmh}(X_i, \, t), \quad t \in [\, 0 \, , 1 \, ].$$

This is a reproducing-kernel regression estimator for pure noise data.

**Lemma 1.** *For every $f \in W^{m,2}(\, 0 \, , 1 \, )$, random or not,*

$$\big| \, \tfrac{1}{n} \sum_{j=1}^{m} D_i \, f(X_i) \, \big| \leq \| \, f \, \|_{wmh} \, \| \, \mathfrak{S}_{nh} \, \|_{wmh},$$

*and under the assumptions (1.4), (1.5) and (1.7), there exists a constant $c_m$ not depending on $h$ such that $\mathbb{E} \big[ \, \| \, \mathfrak{S}_{nh} \, \|_{wmh}^2 \, \big] \leq c_m \, (nh)^{-1}$.*

*Proof.* The identity $\tfrac{1}{n} \sum_{i=1}^{n} D_i \, f(X_i) = \big\langle \, f \, , \, \mathfrak{S}_{nh} \, \big\rangle_{wmh}$ implies the first bound by way of Cauchy-Schwarz. For the expectation, we have

$$\mathbb{E} \big[ \, D_i \, \mathfrak{R}_{wmh}(X_i, \, t) \, \big] = \mathbb{E} \big[ \, \mathbb{E} [\, D_i \mid X_i \,] \, \mathfrak{R}_{wmh}(X_i, \, t) \, \big] = 0$$

and so, since $D_i \, \mathfrak{R}_{wmh}(X_i, \, t), \; i = 1, 2, \ldots, n$, are independent and identically distributed (iid), it follows that

$$\begin{aligned} \mathbb{E} \big[ \, \| \, \mathfrak{S}_{nh} \, \|_{L^2((\, 0 \, , 1 \, ), w)}^2 \, \big] &= n^{-2} \sum_{i=1}^{n} \mathbb{E} \big[ \, D_i^2 \, \| \, \mathfrak{R}_{wmh}(X_i, \, \cdot \,) \, \|_{L^2((\, 0 \, , 1 \, ), w)}^2 \, \big] \\ &\leq n^{-1} \, M \, \mathbb{E} \big[ \, \| \, \mathfrak{R}_{wmh}(X, \, \cdot \,) \, \|_{L^2((\, 0 \, , 1 \, ), w)}^2 \, \big]. \end{aligned}$$

where

$$(2.15) \qquad M = \sup_{x \in [\, 0 \, , 1 \, ]} \, \mathbb{E} [\, D^2 \mid X = x \,].$$

By (1.7), we have $M < \infty$.

Similarly, since $D_i \, \mathfrak{R}_{wmh}^{(m)}(X_i, t), \; i = 1, 2, \ldots, n$, are iid, then

$$\mathbb{E}[\, \| \, \mathfrak{S}_{nh}^{(m)} \, \|_{L^2(0,1)}^2 \,] \leq n^{-1} \, M \, \mathbb{E}[\, \| \, \mathfrak{R}_{wmh}^{(m)}(X, \, \cdot\,) \, \|_{L^2(0,1)}^2 \,].$$

It follows that

$$\mathbb{E}[\, \| \, \mathfrak{S}_{nh} \, \|_{wmh}^2 \,] \leq n^{-1} \, M \, \mathbb{E}[\, \| \, \mathfrak{R}_{wmh}(X, \, \cdot\,) \, \|_{wmh}^2 \,].$$

Now, (2.12) takes care of the last norm. $\qquad\square$

## 3. Random sums

In this section, we discuss sharp bounds on the "random sums" $\mathfrak{S}_{nh}$ of (2.14), using results of Einmahl and Mason [12] regarding convolution-kernel estimators (in a more general setting). Thus, let

$$(3.1) \qquad K \in L^1(\mathbb{R}) \cap L^\infty(\mathbb{R}), \qquad \int_{\mathbb{R}} K(x) \, dx = 1.$$

We also need some restrictions on the "size" of the set of functions on $[\,0\,,1\,]$,

$$(3.2) \qquad \mathcal{K} = \big\{\, K\big(h^{-1}(\,x - \,\cdot\,)\big) \;\big|\; x \in [\,0\,,1\,], \; 0 < h \leq 1 \,\big\}.$$

First, we need to assume that

$$(3.3) \qquad\qquad \mathcal{K} \quad \text{is pointwise measurable}.$$

For the definition of *pointwise measurability*, see van der Vaart and Wellner [23].

Let $Q$ be a probability measure on $([\,0\,,1\,], \mathcal{B})$, and let $\| \, \cdot \, \|_Q$ denote the $L^2(Q)$ metric. For $\varepsilon > 0$, let $\mathcal{N}(\varepsilon, \mathcal{K}, \| \, \cdot \, \|_Q)$ denote the smallest number of balls in the $\| \, \cdot \, \|_Q$ metric needed to cover $\mathcal{K}$, i.e.,

$$(3.4) \qquad \mathcal{N}\big(\varepsilon, \mathcal{K}, \| \cdot \|_Q\big) = \min \left\{\, n \in \mathbb{N} \;\middle|\; \begin{array}{c} \exists \, g_1, g_2, \ldots, g_n \in \mathcal{K} \quad \forall k \in \mathcal{K} \\ \min_{1 \leq i \leq n} \| \, k - g_i \, \|_Q \leq \varepsilon \end{array} \right\}.$$

Then, let

$$(3.5) \qquad \mathcal{N}(\varepsilon, \mathcal{K}) = \sup \mathcal{N}\big(\varepsilon, \mathcal{K}, \| \, \cdot \, \|_Q\big),$$

where the supremum is over *all* probability measures $Q$ on $([\,0\,,1\,], \mathcal{B})$.

The restriction on the size of $\mathcal{K}$ now takes the form that there exist positive constants $C$ and $\nu$ such that

$$(3.6) \qquad \mathcal{N}(\varepsilon, \mathcal{K}) \leq C \, e^{-\nu}, \quad 0 < \varepsilon < 1.$$

Nolan and Pollard [19], see also [23], show that the condition (3.6) holds if the kernel $K$ satisfies (3.1) and (3.3), and has bounded variation,

$$(3.7) \qquad\qquad K \in BV(\mathbb{R}).$$

Whenever $K$ has left and right limits everywhere (so in particular, when (3.7) holds), then (3.3) holds also.

The object of study is the following kernel "estimator" with "pure noise" data,

$$(3.8) \qquad \mathbb{S}_{nh}(t) \stackrel{\text{def}}{=} \frac{1}{n} \sum_{i=1}^{n} D_i \, K_h(X_i - t), \quad t \in [\,0\,,1\,].$$

We quote the following slight modification as it applies to (3.8) of Proposition 2 of Einmahl and Mason [12] without proof. The modification involves the omission of the condition of compact support of the kernel $K$, which is permissible since the design is contained in a compact set, to wit the interval $[\,0\,,1\,]$, [17]. Recall the definition of $\mathcal{H}_n(\gamma)$ from (1.10).

**Proposition 1** (after Einmahl and Mason [12]). *Under the assumptions* (3.1), (3.3), (3.6), (3.7), *and* (1.4), (1.5), *and* (1.7), *for every* $\gamma > 0$,

$$\limsup_{n \to \infty} \sup_{h \in \mathcal{H}_n(\gamma)} \frac{\| \, \mathbb{S}_{nh} \, \|_{\infty}}{\sqrt{(nh)^{-1} \left\{ \log(1/h) \vee \log \log n \right\}}} < \infty \quad almost \ surely.$$

*Proof.* The proof needs updating in only one spot, viz. the bound (3.20) of the Einmahl and Mason [12] paper needs to be established under the present conditions. However, that just amounts to showing that

$$\sup_{0 < h \le 1} \sup_{t \in [\,0\,,1\,]} h \, \mathbb{E}[\,|\, D \, K_h(t - X)\,|^2\,] < \infty.$$

Observe that

$$\mathbb{E}[\,|\, D \, K_h(t - X)\,|^2\,] = \mathbb{E}[\,\mathbb{E}[\, D^2 \mid X\,] \,|\, K_h(t - X)\,|^2\,] \le M \, \mathbb{E}[\,|\, K_h(t - X)\,|^2\,],$$

with $M$ as in (2.15). Now,

$$\mathbb{E}[\,|\, K_h(t - X)\,|^2\,] = \int_0^1 h^{-2} \, K^2\big( h^{-1}(t - x) \big) \, w(x) \, dx$$

$$\le w_2 \, h^{-1} \int_{\mathbb{R}} K^2(x) \, dx \le w_2 \, \| \, K \, \|_{L^1(\mathbb{R})} \, \| \, K \, \|_{L^{\infty}(\mathbb{R})} \, h^{-1} \le c \, h^{-1},$$

for a suitable constant $c$, not depending on $t$. $\qquad \square$

Now, we have the task of relating the random sums involving the reproducing kernels to sums involving convolution kernels. Obviously, some convolution-kernel-like properties of the reproducing kernel are required.

**Definition 1.** We say a family $A_h$, $0 < h < 1$, defined on $[\,0\,,1\,] \times [\,0\,,1\,]$, is convolution-like if it satisfies the following conditions: There exists a constant $c$ such that for all $t \in [\,0\,,1\,]$ and all $h$, $0 < h < 1$,

$$\| \, A_h(\,\cdot\,,t) \, \|_{L^1(0,1)} \le c, \ \| \, A_h(t,\,\cdot\,) \, \|_{\infty} \le c \, h^{-1}, \ |\, A_h(\,\cdot\,,t)\,|_{BV} \le c \, h^{-1}.$$

Here, $|\, f \,|_{BV}$ denotes the total variation of the function $f$ over $[\,0\,,1\,]$.

The families $h^\ell \, \mathfrak{R}_{wmh}^{(\ell)}(t\,,s)$, $0 < h < 1$, $\ell = 0, 1, \ldots, m$, are indeed convolution-like, as shown in [11]. Here,

$$(3.9) \qquad \mathfrak{R}_{wmh}^{(\ell)}(t\,,s) = \frac{d^\ell}{ds^\ell} \, \mathfrak{R}_{wmh}(t\,,s)$$

denotes the $\ell$-th order derivative of $\mathfrak{R}_{wmh}(t, s)$ with respect to $s$ (or by symmetry, with respect to $t$). This result is in the style of results on the "equivalent" kernel for spline smoothing, except that that the kernel is *not* a convolution kernel and that it handles arbitrary design densities subject to the condition (1.5) *and* treats the boundary conditions in (2.13) exactly. The relevant references on equivalent kernels for spline smoothing are [1, 5, 6, 7, 11, 18, 20, 21].

Now, there is an interesting way of connecting the reproducing kernel sum $\mathfrak{S}_{nh}$ to a sum $\mathbb{S}_{nh}$ for an appropriate kernel $K$. Define

$$(3.10) \qquad \begin{aligned} g(x) &= \exp(-x)\, \mathbb{1}(x \geq 0), \\ g_h(x) &= h^{-1} g(h^{-1}x), \quad x \in \mathbb{R}. \end{aligned}$$

One verifies that $h\, g_h$ is the fundamental solution for the initial value problem

$$(3.11) \qquad \begin{aligned} h\, u' + u &= v \quad \text{on } (0, 1), \\ u(0) &= a, \end{aligned}$$

i.e., for $1 \leq p \leq \infty$ and $v \in L^p(0, 1)$, the solution $u$ of the initial value problem (3.11) satisfies $u \in L^p(0, 1)$, and is given by

$$(3.12) \qquad u(x) = h\, g_h(x)\, u(0) + \int_0^1 g_h(x - z)\, v(z)\, dz, \quad x \in [0, 1],$$

see, e.g., [4], Section 2.1, formulas (10) through (14). Note that the last integral is really only over the interval $[0, x]$. Since $v = h\, u' + u$, this leads to the integral representation of the function $u$,

$$(3.13) \qquad u(x) = h\, g_h(x)\, u(0) + \int_0^1 g_h(x - z)\, \{h\, u'(z) + u(z)\}\, dz, \quad x \in [0, 1].$$

Now, one verifies that

$$(3.14) \qquad \text{the kernel } g \text{ satisfies (3.1), (3.3), (3.7),}$$

so that the class $\Gamma$ generated by $g$,

$$(3.15) \qquad \Gamma = \{g(h^{-1}(x - \cdot)) \mid x \in [0, 1], \, 0 < h \leq 1\}.$$

satisfies (3.6), i.e.,

$$(3.16) \qquad \mathcal{N}(\varepsilon, \Gamma) \leq C\, e^{-\nu}, \quad 0 < \varepsilon < 1.$$

Thus, Proposition 1 would apply to the random sum

$$(3.17) \qquad s^{nh}(z) = \frac{1}{n} \sum_{i=1}^{n} D_i\, g_h(X_i - z), \quad z \in [0, 1],$$

but first we connect the sums $\mathfrak{S}_{nh}$ and $s^{nh}$.

**Lemma 2.** *Assume that the functions $A_h$, $0 < h \leq 1$, are convolution-like in the sense of Definition 1. Then, there exists a constant $c$ such that for all $h$, $0 < h \leq 1$, all $D_1, D_2, \ldots, D_n \in \mathbb{R}$, and all positive $X_1, X_2, \ldots, X_n \in [0, 1]$,*

$$\left\| \frac{1}{n} \sum_{i=1}^{n} D_i\, A_h(X_i, \cdot) \right\|_{\infty} \leq c \left\| \frac{1}{n} \sum_{i=1}^{n} D_i\, g_h(X_i - \cdot) \right\|_{\infty}.$$

*Proof.* Assume that $A_h$ is differentiable with respect to its first argument. Then,

$$| A_h(\,\cdot\,, t)\,|_{BV} = \|\, A_h'(\,\cdot\,, t)\,\|_{L^1(0,1)},$$

where the prime $\prime$ denotes differentiation with respect to the first argument. Now, apply (3.13) to the function $u = A_h(\,\cdot\,, t)$ (for fixed $t$), so for all $x$,

$$(3.18) \quad A_h(x, t) = h\, g_h(x)\, A_h(0, t) + \int_0^1 g_h(x-z)\, \left\{\, h\, A_h'(z, t) + A_h(z, t)\,\right\} dz.$$

Next, take $x = X_i$ and substitute this into

$$S_{nh}(t) \overset{\text{def}}{=} \tfrac{1}{n} \sum_{i=1}^n D_i\, A_h(X_i, t), \quad t \in [\,0, 1\,].$$

Then, we have

$$(3.19) \qquad S_{nh}(t) = h\, A_h(0, t)\, s^{nh}(0) + \int_0^1 \left\{\, h\, A_h'(z, t) + A_h(z, t)\,\right\} s^{nh}(z)\, dz.$$

Now, straightforward bounding gives

$$\|\, S_{nh}\,\|_\infty \le C \,|\, s^{nh}(0)\,| + C_1 \,\|\, s^{nh}\,\|_\infty,$$

where $C = h\, \|\, A_h(0, \,\cdot\,)\,\|_\infty$ and

$$C_1 = \sup_{t \in [\,0,1\,]} \|\, A_h(t, \,\cdot\,)\,\|_{L^1(0,1)} + h\, |\, A_h(t, \,\cdot\,)\,|_{BV}.$$

So, by the convolution-like properties of $A_h$, the constants $C$ and $C_1$ are bounded, uniformly in $h$. Also, since all of the $X_i$ are positive, then

$$\lim_{z \to 0} s^{nh}(z) = s^{nh}(0),$$

and so, for $C_2 = C + C_1$, we have $\|\, S_{nh}\,\|_\infty \le C_2 \,\|\, s^{nh}\,\|_\infty$.

The extension to the case where $A_h$ is not necessarily differentiable with respect to its first argument follows readily. $\qquad \square$

Since the families $h^\ell\, \mathfrak{R}^{(\ell)}_{wmh}(t, s)$, $0 < h < 1$, $\ell = 0, 1, \ldots, m$, are convolution-like in the sense of Definition 1, we may apply the lemma to the sum $\mathfrak{S}_{nh}$ of (2.14) and its derivatives. This yields

$$(3.20) \qquad h^\ell\, \|\, \mathfrak{S}^{(\ell)}_{nh}\,\|_\infty \le c\, \|\, s^{nh}\,\|_\infty, \quad \ell = 0, 1, \ldots, m.$$

Now, for the model (1.1) through (1.7), the sum $s^{nh}$ of (3.17) may be treated by the above formulated Proposition 1. This proves the following result. (Recall the definition (1.10) of $\mathcal{H}_n(\gamma)$.)

**Theorem 2.** *Under the assumptions* (1.4), (1.5) *and* (1.7), *for* $\gamma > 0$, *and for* $\ell = 0, 1, \ldots, m$,

$$Q_{\infty,\ell}(\gamma) \overset{\text{def}}{=} \limsup_{n \to \infty}\, \sup_{h \in \mathcal{H}_n(\gamma)}\, \frac{h^\ell\, \|\, \mathfrak{S}^{(\ell)}_{nh}\,\|_\infty}{\sqrt{(nh)^{-1}\,\{\log(1/h) \vee \log\log n\}}} < \infty,$$

*almost surely.*

It turns out that we need a similar result for the $\| \cdot \|_{wmh}$ norm, which requires a result for the $L^2$ norm. The following is good enough for our purposes. Obviously, with $\mathfrak{S}_{nh}^{(m)}$ denoting the $m$-th order derivative of $\mathfrak{S}_{nh}$, we have

$$\| \mathfrak{S}_{nh} \|_{L^2((0,1),w)} \le c \, \| \mathfrak{S}_{nh} \|_{\infty}, \quad \| \mathfrak{S}_{nh}^{(m)} \|_{L^2(0,1)} \le \| \mathfrak{S}_{nh}^{(m)} \|_{\infty},$$

with $c = \sqrt{w_2}$, and then Theorem 2 gives useful bounds for the $\| \cdot \|_{wmh}$ norm.

**Corollary 3.** *Under the conditions of Theorem 2, we have almost surely,*

$$Q_{wm} \overset{\text{def}}{=} \limsup_{n \to \infty} \sup_{h \in \mathcal{H}_n(\gamma)} \frac{\| \mathfrak{S}_{nh} \|_{wmh}}{\sqrt{(nh)^{-1} \{ \log(1/h) \vee \log \log n \}}} < \infty.$$

## 4. The design sums

The reproducing kernel Hilbert space set-up is also useful for connecting *random design sums*

$$\tfrac{1}{n} \sum_{i=1}^n | f(X_i) |^2$$

to their (partial?) expectations

$$\int_0^1 | f(x) |^2 \, w(x) \, dx,$$

for *random* functions $f$. In particular, we prove the following almost sure result. The range of the smoothing parameter is much larger here than it was before, although we only need it for $h \in \mathcal{H}_n(\gamma)$, see (1.10). Here, let

(4.1) $$\mathcal{D}_n(\gamma) = \left[ \, \gamma \, n^{-1} \log n \, , \, \tfrac{1}{2} \, \right].$$

**Theorem 4.** *Under the assumptions* (1.4) *and* (1.5), *for all* $f \in W^{m,2}(0,1)$,

$$\tfrac{1}{n} \sum_{j=1}^n | f(X_i) |^2 + h^{2m} \| f^{(m)} \|^2 \ge r_{nh} \| f \|_{wmh}^2,$$

*where* $\quad \liminf_{n \to \infty} \inf_{h \in \mathcal{D}_n(\gamma)} r_{nh} = 1 \quad$ *almost surely.*

To prove this, let $W$ be the (cumulative) distribution function corresponding to the design density $w$ and let $W_n$ be the empirical distribution function of the design $X_1, X_2, \ldots, X_n$, and introduce the "design sums"

(4.2) $$w^{nh}(t) \overset{\text{def}}{=} g_h * dW_n(t) \overset{\text{def}}{=} \tfrac{1}{n} \sum_{i=1}^n g_h(X_i - t), \; t \in [0,1],$$

which is a *convolution*-kernel density estimator and its expectation,

(4.3) $$\mathbb{E}[\, w^{nh}(t) \,] = g_h * dW(t) = \int_0^1 g_h(\tau - t) \, w(\tau) \, d\tau, \quad t \in [0,1].$$

We will use Theorem 1 of Einmahl and Mason [12], quoted here for convenience. (This time, no modifications are necessary.)

**Proposition 2** ([12]). *Under the assumptions* (3.1), (3.3), (3.6), *and* (3.7), *and* (1.4) *and* (1.5), *for every* $\gamma > 0$,

$$\limsup_{n \to \infty} \sup_{h \in \mathcal{D}_n(\gamma)} \frac{\| w^{nh} - \mathbb{E}[w^{nh}] \|_\infty}{\sqrt{nh \left( \log(1/h) \vee \log \log n \right)}} < \infty \quad \text{almost surely.}$$

To prove Theorem 4, we start with simple "design sums".

**Lemma 3.** *Under the assumptions* (1.4) *and* (1.5), *for all* $f \in W^{1,1}(0, 1)$,

$$\left| \int_0^1 f(t) \left\{ dW_n(t) - dW(t) \right\} \right| \le \zeta_{nh} \left\{ \| f \|_{L^1((0,1),w)} + h \| f' \|_{L^1(0,1)} \right\}$$

*where*

$$\limsup_{n \to \infty} \sup_{h \in \mathcal{D}_n(\gamma)} \frac{\zeta_{nh}}{\sqrt{(nh)^{-1} \left\{ \log(1/h) \vee \log \log n \right\}}} < \infty$$

*almost surely.*

*Proof.* With the reproducing kernel Hilbert space trick,

$$f(t) = \left\langle f, \, \mathfrak{R}_{w,1,h}(t, \cdot) \right\rangle_{w,1,h},$$

we obtain by linearity and Fubini's theorem that

$$\int_0^1 f(t) \left\{ dW_n(t) - dW(t) \right\} = \left\langle f, \delta^{nh} \right\rangle_{w,1,h},$$

where

$$\delta^{nh}(s) = \int_0^1 \mathfrak{R}_{w,1,h}(t, s) \left\{ dW_n(t) - dW(t) \right\}, \quad s \in [0, 1],$$

is the variance part of the pointwise error of a reproducing-kernel estimator of the design density $w$. Now, straightforward bounding gives

$$\left\langle f, \delta^{nh} \right\rangle_{L^2((0,1),w)} \le \| f \|_{L^1((0,1),w)} \| \delta^{nh} \|_\infty,$$

$$\left\langle f', (\delta^{nh})' \right\rangle_{L^2(0,1)} \le \| f' \|_{L^1(0,1)} \| (\delta^{nh})' \|_\infty,$$

so that

$$\left\langle f, \delta^{nh} \right\rangle_{w,1,h} \le \left\{ \| f \|_{L^1((0,1),w)} + h \| f' \|_{L^1(0,1)} \right\} \left\{ \| \delta^{nh} \|_\infty + h \| (\delta^{nh})' \|_\infty \right\},$$

with, explicitly,

$$\| \delta^{nh} \|_\infty = \left\| \tfrac{1}{n} \sum_{i=1}^n \mathfrak{R}_{w,1,h}(X_i, \cdot) - \mathbb{E}[\mathfrak{R}_{w,1,h}(X_1, \cdot)] \right\|_\infty,$$

$$h \| (\delta^{nh})' \|_\infty = \left\| \tfrac{1}{n} \sum_{i=1}^n h \, \mathfrak{R}'_{w,1,h}(X_i, \cdot) - h \, \mathbb{E}[\mathfrak{R}'_{w,1,h}(X_1, \cdot)] \right\|_\infty.$$

Both of these may be interpreted as the variance part of the uniform error of (reproducing) kernel estimators. As already noted, the families $\mathfrak{R}_{wmh}(t, s)$ and $h \mathfrak{R}'_{w,1,h}(t, s)$ are convolution-like in the sense of Definition 1. Then, by an appeal to Lemma 2 with $D_i = 1$ for all $i$,

$$\| \delta^{nh} \|_\infty \le C \| g_h * \left\{ dW_n - dW \right\} \|_\infty,$$

$$\| (\delta^{nh})' \|_\infty \le C_1 \| g_h * \left\{ dW_n - dW \right\} \|_\infty,$$

for suitable constants $C$ and $C_1$. Now, an appeal to Theorem 1 of Einmahl and Mason [12], see Proposition 2 above, clinches the deal. $\qquad\square$

We now get the following lemma, which immediately implies Theorem 4.

**Lemma 4.** *Under the assumptions* (1.4) *and* (1.5), *for all* $f$, $g \in W^{m,2}(0,1)$,

$$\Big| \int_0^1 f(t)\, g(t) \left\{ dW_n(t) - dW(t) \right\} \Big| \leq \eta_{nh} \, \|f\|_{wmh} \, \|g\|_{wmh},$$

*where* $\quad \limsup_{n \to \infty} \; \sup_{h \in \mathcal{D}_n(\gamma)} \; \dfrac{\eta_{nh}}{\sqrt{(nh)^{-1} \left\{ \log(1/h) \vee \log \log n \right\}}} < \infty$

*almost surely.*

*Proof.* From Lemma 3, we get the bound

$$\zeta_{nh} \left\{ \| f\,g \|_{L^1((0,1),w)} + h \, \| (f\,g)' \|_{L^1(0,1)} \right\},$$

with the requisite behavior of $\zeta^{nh}$. Now, from (2.8),

$$\| f\,g \|_{L^1((0,1),w)} + h \, \| (f\,g)' \|_{L^1(0,1)} \leq c \, \|f\|_{w,1,h} \, \|g\|_{w,1,h}$$

for an appropriate constant $c$. Finally, (2.9) gives $\|f\|_{w,1,h} \leq c \, \|f\|_{w,m,h}$, again for an appropriate constant $c$, and likewise for $g$. Thus, $\eta_{nh}$ satisfies $\eta_{nh} \leq c \, \zeta_{nh}$, and the lemma follows. $\qquad\square$

## 5. $L^2$ error bounds

We are now ready to prove almost sure bounds on $\| f^{nh} - f_o \|_{wmh}^2$ for the spline smoother $f^{nh}$. The starting point is the quadratic Taylor expansion of the objective function $\mathrm{IS}(f)$ of (1.8) around its minimizer. Let

$$(5.1) \qquad\qquad \varepsilon \equiv f^{nh} - f_o.$$

Since the Gateaux variation of IS at its minimizer vanishes, this gives

$$(5.2) \qquad \frac{1}{n} \sum_{i=1}^n |\varepsilon(X_i)|^2 + h^{2m} \|\varepsilon^{(m)}\|^2 = \mathrm{IS}(f_o) - \mathrm{IS}(f^{nh}).$$

Now, again, simple quadratic Taylor expansion around $f_o$ gives

$$(5.3) \qquad \begin{aligned} \mathrm{IS}(f_o) - \mathrm{IS}(f^{nh}) = &-\frac{1}{n} \sum_{i=1}^n |\varepsilon(X_i)|^2 + \frac{2}{n} \sum_{i=1}^n D_i\, \varepsilon(X_i) \\ &- h^{2m} \|\varepsilon^{(m)}\|^2 + 2\,h^{2m} \left\langle f_o^{(m)}, \, \varepsilon^{(m)} \right\rangle, \end{aligned}$$

and so, after substitution into (5.2),

$$(5.4) \qquad \frac{1}{n} \sum_{i=1}^n |\varepsilon(X_i)|^2 + h^{2m} \|\varepsilon^{(m)}\|^2 = \frac{1}{n} \sum_{i=1}^n D_i\, \varepsilon(X_i) + h^{2m} \left\langle f_o^{(m)}, \, \varepsilon^{(m)} \right\rangle.$$

This is similar to the development in [22].

Now, with Lemma 1, Theorem 4, and Cauchy-Schwarz, one obtains

$$(5.5) \qquad r^{nh} \, \| \, \varepsilon \, \|^2_{wmh} \leq \| \, \varepsilon \, \|_{wmh} \, \big\{ \, \| \, \mathfrak{S}_{nh} \, \|_{wmh} + h^m \, \| \, f_o^{(m)} \, \| \, \big\},$$

where we took the liberty of using $h^m \| \, \varepsilon^{(m)} \, \| \leq \| \, \varepsilon \, \|_{wmh}$. It follows that

$$(5.6) \qquad r^{nh} \, \| \, \varepsilon \, \|_{wmh} \leq \| \, \mathfrak{S}_{nh} \, \|_{wmh} + h^m \, \| \, f_o^{(m)} \, \|,$$

and the following result emerges.

**Theorem 5.** *For the model* (1.1), *under the assumptions* (1.4), (1.6), (1.2), *and* (1.5), *with* $\mathcal{H}_n(\gamma)$ *defined in* (1.10), *for* $\gamma > 0$, *almost surely,*

$$\limsup_{n \to \infty} \; \sup_{h \in \mathcal{H}_n(\gamma)} \; \frac{\| \, f^{nh} - f_o \, \|_{wmh}}{\sqrt{h^{2m} + (nh)^{-1} \, \{ \, \log(1/h) \vee \log \log n \, \}}} < \infty,$$

*and for* $h \asymp (n^{-1} \log n)^{1/(2m+1)}$ *(deterministic or random),*

$$\| \, f^{nh} - f_o \, \|_{wmh} = \mathcal{O}\big( \, (n^{-1} \log n)^{m/(2m+1)} \, \big) \quad \text{almost surely.}$$

*Proof.* This follows from (5.6) and Corollary 3. $\qquad \square$

The error bound (5.6) appears to be quite sharp. In the next section, we show that $f^{nh}(t) - \mathbb{E}[\, f^{nh}(t) \mid X_1, \ldots, X_n \,] \approx \mathfrak{S}_{nh}$ in a precise, useful sense.

## 6. C-splines

In this section, we determine a useful, accurate expression for the variance part $f^{nh}(t) - \mathbb{E}[\, f^{nh}(t) \mid X_1, \ldots, X_n \,]$ of the pointwise error $f^{nh}(t) - f_o(t)$, with an eye towards almost sure uniform error bounds. Since the estimator $f^{nh}$ is linear in the data, one sees that

$$(6.1) \qquad \varphi^{nh} = f^{nh} - \mathbb{E}[\, f^{nh} \mid X_1, \ldots, X_n \,]$$

is the solution to the "pure noise" problem

$$(6.2) \qquad \begin{aligned} &\text{minimize} \quad \frac{1}{n} \sum_{i=1}^{n} | \, f(X_i) - D_i \, |^2 + h^{2m} \, \| \, f^{(m)} \, \|^2 \\ &\text{subject to} \quad f \in W^{m,2}(\, 0\,,\, 1\,). \end{aligned}$$

In fact, we show that $\varphi^{nh}(t) \approx \psi^{nh}(t)$, where $f = \psi^{nh}$ solves the C(ontinuous)-spline problem

$$(6.3) \qquad \begin{aligned} &\text{minimize} \quad \| \, f \, \|^2_{L^2((\,0\,,\,1\,), w)} - \frac{2}{n} \sum_{i=1}^{n} D_i \, f(X_i) + h^{2m} \, \| \, f^{(m)} \, \|^2 \\ &\text{subject to} \quad f \in W^{m,2}(\, 0\,,\, 1\,). \end{aligned}$$

By the interpretation of the reproducing kernel $\mathfrak{R}_{wmh}$ as the Green's function for the boundary value problem (2.13), one observe that $\psi^{nh}$ is given by

$$(6.4) \qquad \psi^{nh}(t) = \frac{1}{n} \sum_{i=1}^{n} D_i \, \mathfrak{R}_{wmh}(X_i, \, t),$$

so that the following almost sure error bounds apply.

**Theorem 6.** *Under the assumptions of Theorem 5, almost surely, uniformly in* $h \in H_n(\gamma)$, *see* (1.10),

$$\| \varphi^{nh} - \psi^{nh} \|_{wmh} = \mathcal{O}\big( h^{-1/2} (nh)^{-1} \{ \log(1/h) \vee \log\log n \} \big),$$
$$\| \varphi^{nh} - \psi^{nh} \|_{\infty} = \mathcal{O}\big( h^{-1} (nh)^{-1} \{ \log(1/h) \vee \log\log n \} \big).$$

*Proof.* Let $\varepsilon = \varphi^{nh} - \psi^{nh}$. Similar to the inequality (5.4), one obtains quadratic inequalities for the discrete and the continuous spline problems. Adding these gives

$$(6.5) \qquad \| \varepsilon \|^2 + 2\, h^{2m} \| \varepsilon^{(m)} \|^2 + \tfrac{1}{n} \sum_{i=1}^n | \varepsilon(X_i) |^2 = \text{rhs},$$

where

$$\text{rhs} = \int_0^1 g(t) \{ dW_n(t) - dW(t) \},$$

with $g = | \varphi^{nh} |^2 - | \psi^{nh} |^2 = ( \varphi^{nh} + \psi^{nh} ) \varepsilon$, and $\varepsilon$ as above.

Now let $\gamma > 0$ be fixed. Then, the following statements hold uniformly in $h \in \mathcal{H}_n(\gamma)$. Using Lemma 3, one obtains, almost surely

$$\text{rhs} = \mathcal{O}\big( \sqrt{(nh)^{-1} \{ \log(1/h) \vee \log\log n \}} \, \big) \{ \| g \|_{L^1((0,1),w)} + h \| g \|_{L^1(0,1)} \}$$
$$= \mathcal{O}\big( \sqrt{(nh)^{-1} \{ \log(1/h) \vee \log\log n \}} \, \big) \| \varepsilon \|_{wmh} \| \varphi^{nh} + \psi^{nh} \|_{wmh},$$

where we used the multiplication result (2.8), and (2.9). Substituting this into (6.5), one obtains almost surely,

$$\| \varepsilon \|_{wmh} = \mathcal{O}\big( \sqrt{(nh)^{-1} \{ \log(1/h) \vee \log\log n \}} \, \big) \| \varphi^{nh} + \psi^{nh} \|_{wmh}.$$

Now,

$$\| \varphi^{nh} + \psi^{nh} \|_{wmh} \leq \| \varphi^{nh} \|_{wmh} + \| \psi^{nh} \|_{wmh},$$

and consequently, by Theorem 5 applied with $f_o = 0$, and (6.5),

$$\| \varphi^{nh} + \psi^{nh} \|_{wmh} = \mathcal{O}\big( \sqrt{(nh)^{-1} \{ \log(1/h) \vee \log\log n \}} \, \big) \quad \text{almost surely}.$$

Thus,

$$\| \varepsilon \|_{wmh} = \mathcal{O}\big( (nh)^{-1} \{ \log(1/h) \vee \log\log n \} \big),$$

and so, at the loss of a factor $h^{-1/2}$,

$$\| \varepsilon \|_{\infty} = \mathcal{O}\big( h^{-1/2}(nh)^{-1} \{ \log(1/h) \vee \log\log n \} \big).$$

The theorem has been proved. ☐

The above completes the proof of Theorem 1.

## 7. C-splines: the Bias

In considering the bias of the estimator $f^{nh}$, note that $f_h = \mathbb{E}[\, f^{nh} \mid X_1, \ldots, X_n]$ is the solution of (1.8) with $D_i = 0$, $i = 1, 2, \ldots, n$, i.e., the solution of the discrete noiseless problem

$$(7.1) \qquad \text{minimize} \quad DN(f) \stackrel{\text{def}}{=} \tfrac{1}{n} \sum_{i=1}^n | f(X_i) - f_o(X_i) |^2 + h^{2m} \| f^{(m)} \|^2$$
$$\text{subject to} \quad f \in W^{m,2}(0,1).$$

Note that the randomness in $f_h$ is due to the randomness of the design. We must compare $f_h$ to $f_o$, but it is easier to first compare $f_h$ to $\varphi_h$, the solution of the continuous noiseless problem

$$(7.2) \qquad \begin{aligned} &\text{minimize} \quad CN(f) \stackrel{\text{def}}{=} \|f - f_o\|^2_{L^2((0,1),w)} + h^{2m} \|f^{(m)}\|^2 \\ &\text{subject to} \quad f \in W^{m,2}(0,1). \end{aligned}$$

In this section we prove the following theorem on the conditional bias. Note the "restricted" set $\mathcal{G}_n(\gamma)$ of allowable values of $h$.

**Theorem 7.** *Under the assumptions of Theorem 5, almost surely,*

$$\limsup_{n \to \infty} \sup_{h \in \mathcal{G}_n(\gamma)} \frac{\|f_h - f_o\|_\infty}{\sqrt{h^{2m} + (nh)^{-1}\{\log(1/h) \vee \log \log n\}}} < \infty.$$

The proof goes again by way of the reproducing kernel approximation, and follows without further ado from the following two lemmas.

**Lemma 5.** *Under the assumptions of Theorem 5, almost surely,*

$$\limsup_{n \to \infty} \sup_{h \in \mathcal{H}_n(\gamma)} \frac{\|\varphi_h - f_h\|_\infty}{h^{-1/2}\{h^{2m} + (nh)^{-1}\{\log(1/h) \vee \log \log n\}\}} < \infty.$$

*Proof.* Let $\varepsilon = \varphi_h - f_h$. Similar to the derivation of (5.6), one obtains

$$(7.3) \qquad\qquad 2\|\varepsilon\|^2_{wmh} = \text{rhs}$$

where "rhs" $= DN(\varphi_h) - DN(f_h) + CN(f_h) - CN(\varphi_h)$. This simplifies to

$$\text{rhs} = \int_0^1 g(t)\{dW_n(t) - dW(t)\},$$

with $W$ and $W_n$ as in Lemma 3, and

$$\begin{aligned} g(t) &= |\varphi_h(t) - f_o(t)|^2 - |f_h(t) - f_o(t)|^2 \\ &= (\varphi_h(t) + f_h(t) - 2f_o(t))\,\varepsilon(t). \end{aligned}$$

By Lemma 4, we get that

$$(7.4) \qquad\qquad \text{rhs} \le \eta_{nh} \|\varphi_h + f_h - 2f_o\|_{wmh} \|\varepsilon\|_{wmh},$$

with

$$(7.5) \qquad\qquad \eta_{nh} = \mathcal{O}\left(\sqrt{(nh)^{-1}\{\log(1/h) \vee \log \log n\}}\right)$$

almost surely, uniformly in $h \in \mathcal{H}_n(\gamma)$.

Substituting (7.4) into (7.3), we obtain

$$(7.6) \qquad\qquad \|\varepsilon\|_{wmh} \le \tfrac{1}{2}\eta_{nh}\|\varphi_h + f_h - 2f_o\|_{wmh}.$$

Now, with regards to bounding $\|\varphi_h + f_h - 2f_o\|_{wmh}$, the situation is as in Section 5, except that the $D_i = 0$, $i = 1, 2, \ldots, n$. Then, almost surely,

$$\limsup_{n \to \infty} \sup_{h \in \mathcal{H}_n(\gamma)} h^{-m} \|f_h - f_o\|_{wmh} < \infty.$$

(Note that there is still randomness in $f_h$ due to the design.) In the same way, one obtains that deterministically,

$$\| \, \varphi_h - f_o \, \|_{wmh} = \mathcal{O}\big( \, h^m \, \big).$$

It then follows from (7.6) that

$$\limsup_{n \to \infty} \ \sup_{h \in \mathcal{H}_n(\gamma)} \ \frac{\| \, \varepsilon \, \|_{wmh}}{h^m \ \sqrt{ (nh)^{-1} \{ \, \log(1/h) \vee \log \log n \, \} }} < \infty.$$

Of course,

$$2 \, h^m \, \sqrt{ (nh)^{-1} \{ \, \log(1/h) \vee \log \log n \, \} } \le h^{2m} + (nh)^{-1} \{ \, \log(1/h) \vee \log \log n \, \},$$

so that

$$\frac{2 \, \| \, \varepsilon \, \|_{wmh}}{h^{2m} + (nh)^{-1} \{ \, \log(1/h) \vee \log \log n \, \}} \le \frac{\| \, \varepsilon \, \|_{wmh}}{h^m \ \sqrt{ (nh)^{-1} \{ \, \log(1/h) \vee \log \log n \, \} }}.$$

At the loss of a factor $h^{1/2}$, this gives us the required bound on $\| \, \varepsilon \, \|_\infty$. □

**Lemma 6.** *Under the assumption* (1.2)*, there exists a constant* $c$*, such that*

$$\| \, \varphi_h - f_o \, \|_\infty \le c \, h^m \, \| \, f_o^{(m)} \, \|_\infty,$$

*provided* $f_o$ *satisfies* (1.5)*.*

*Proof.* One verifies that $\varphi_h$ is the solution to the differential equation

$$(7.7) \qquad\qquad (-h^2)^m \, f^{(2m)} + w \, f = w \, f_o \quad \text{on } ( \, 0 \, , 1 \, ),$$

supplemented with the natural boundary conditions. Now, we assume that the regression function $f_o$ satisfies $f_o \in W^{m,\infty}( \, 0 \, , 1 \, )$, so certainly, $f_o \in W^{m,2}( \, 0 \, , 1 \, )$. Then, cf. (2.13), the solution of (7.7) is given by

$$\varphi_h(t) = \int_0^1 \mathfrak{R}_{wmh}( \, t \, , s \, ) \, w( \, s \, ) \, f_o( \, s \, ) \, ds = \big\langle \, \mathfrak{R}_{wmh}( \, t \, , \cdot \, ) \, , f_o \, \big\rangle_{L^2(( \, 0 \, , 1 \, ), w)},$$

and so

$$\varphi_h(t) = \big\langle \, \mathfrak{R}_{wmh}( \, t \, , \cdot \, ) \, , f_o \, \big\rangle_{wmh} - h^{2m} \big\langle \, \mathfrak{R}_{wmh}^{(m)}( \, t \, , \cdot \, ) \, , f_o^{(m)} \, \big\rangle_{L^2(0,1)},$$

with $\mathfrak{R}_{wmh}^{(m)}( \, t \, , s \, )$ the $m$-th derivative of $\mathfrak{R}_{wmh}( \, t \, , s \, )$ with respect to $s$, as in (3.9). Since, $\big\langle \, \mathfrak{R}_{wmh}( \, t \, , \cdot \, ) \, , f_o \, \big\rangle_{wmh} = f_o(t)$, and

$$\big| \, \big\langle \, \mathfrak{R}_{wmh}^{(m)}( \, t \, , \cdot \, ) \, , f_o^{(m)} \, \big\rangle_{L^2(0,1)} \, \big| \le \| \, \mathfrak{R}_{wmh}^{(m)}( \, t \, , \cdot \, ) \, \|_{L^1(( \, 0 \, , 1 \, ), w)} \, \| \, f_o^{(m)} \, \|_\infty$$
$$\le c \, h^{-m} \, \| \, f_o^{(m)} \, \|_\infty,$$

the last inequality by the convolution-likeness of $h^m \, \mathfrak{R}_{wmh}^{(m)}$, the lemma follows. □

## 8. Uniform error bounds

As an application of the reproducing-kernel approximation to the spline smoother, we obtain uniform error bounds on the spline smoother, uniformly in the bandwidth over a wide (useful) range.

*Proof of the Main Theorem.* First, let us consider the result of Theorem 1. Recall that the range of the bandwidths is $\mathcal{G}_n(\gamma) = \left[ \gamma \, (n^{-1} \log n)^{1-2/\lambda}, \, 1/2 \right]$ and that $2 < \lambda < \min(\kappa, 4)$ with $\kappa > 3$. So, $h \in \mathcal{G}_n(\gamma)$ implies that $h \gg (n^{-1} \log n)^{1/2}$. Now, from Theorem 6, with $f_h = \mathbb{E}[\, f^{nh} \mid X_1, \cdots, X_n\,]$, and $\mathfrak{S}_{nh}$ given by (2.14),

$$f^{nh}(t) - f_h(t) = \mathfrak{S}_{nh}(t) + \varepsilon^{nh}(t),$$

with, almost surely, uniformly in $h \in \mathcal{H}_n(\gamma)$,

$$\| \, \varepsilon^{nh} \, \|_\infty = \mathcal{O}\Big( h^{-1/2}(nh)^{-1} \left\{ \log(1/h) \vee \log \log n \right\} \Big).$$

For $h \gg (n^{-1} \log n)^{1/2}$, we may conclude that

$$\| \, \varepsilon^{nh} \, \|_\infty = o\Big( \sqrt{(nh)^{-1} \left\{ \log(1/h) \vee \log \log n \right\}} \Big),$$

which is negligible compared to the upperbound of Theorem 2,

$$\| \, \mathfrak{S}_{nh} \, \|_\infty = \mathcal{O}\big( \sqrt{(nh)^{-1} \left\{ \log(1/h) \vee \log \log n \right\}} \big)$$

almost surely, uniformly in $h \in \mathcal{H}_n(\gamma)$. Finally,

$$\begin{aligned}
\| \, f^{nh} - f_o \, \|_\infty &\le \| \, f^{nh} - f_h \, \|_\infty + \| \, f_h - f_o \, \|_\infty \\
&\le \| \, \mathfrak{S}_{nh} \, \|_\infty + \| \, \varepsilon^{nh} \, \|_\infty + \| \, f_h - f_o \, \|_\infty,
\end{aligned}$$

and Theorem 7 takes care of the last term.                               □

**Acknowledgment.**   We thank David Mason for patiently explaining the results of Einmahl and Mason [12] to us, and for straightening out the required modification of their Proposition 2.